\documentclass[mathpazo]{cicp}

\usepackage{listings}
\usepackage{xcolor}

\definecolor{codegreen}{rgb}{0,0.6,0}
\definecolor{codegray}{rgb}{0.5,0.5,0.5}
\definecolor{codepurple}{rgb}{0.58,0,0.82}
\definecolor{backcolour}{rgb}{0.95,0.95,0.92}

\lstdefinestyle{mystyle}{
backgroundcolor=\color{backcolour},   
commentstyle=\color{codegreen},
keywordstyle=\color{magenta},
numberstyle=\tiny\color{codegray},
stringstyle=\color{codepurple},
basicstyle=\ttfamily\footnotesize,
breakatwhitespace=false,         
breaklines=true,                 
captionpos=b,                    
keepspaces=true,                 
numbers=left,                    
numbersep=5pt,                  
showspaces=false,                
showstringspaces=false,
showtabs=false,                  
tabsize=2
}
\usepackage{amsmath,epsf,cite}
\usepackage{amssymb,amsthm,tocvsec2}
\usepackage{bigints}
\usepackage{graphicx}
\usepackage{sidecap}
\usepackage{wrapfig}
\usepackage{pict2e}
\usepackage{float}
\usepackage{graphicx}
\usepackage{setspace}
\usepackage{float}
\usepackage{color}
\usepackage{mathrsfs}
\usepackage{epsfig,epsf,latexsym,subfigure}
\usepackage{fancybox}
\usepackage{fancyhdr}
\usepackage{mdwlist}
\usepackage{paralist}
\usepackage{setspace}
\usepackage{tikz}
\usetikzlibrary{shapes,arrows}
\usepackage{enumitem}
\setlist{nolistsep}
\usepackage{indentfirst}
\usepackage{epstopdf}

\numberwithin{equation}{section}
\numberwithin{figure}{section}
\numberwithin{table}{section}

\theoremstyle{plain}

\theoremstyle{remark}

\newcommand{\ben}{\begin{eqnarray}}
\newcommand{\een}{\end{eqnarray}}
\newcommand{\benx}{\begin{eqnarray*}}
\newcommand{\eenx}{\end{eqnarray*}}
\newcommand{\beq}{\begin{equation}}
\newcommand{\eeq}{\end{equation}}
\newcommand{\beqx}{\begin{equation*}}
\newcommand{\eeqx}{\end{equation*}}
\newcommand{\bea}{\begin{array}}
\newcommand{\eea}{\end{array}}
\newcommand{\bef}{\begin{figure}[H]}
\newcommand{\eef}{\end{figure}}
\newcommand{\be}{\begin{equation}}
\newcommand{\ee}{\end{equation}}

\newcommand{\bse}{\begin{subequations}}
\newcommand{\ese}{\end{subequations}}

\def\hR{\mathbb R}
\def\bX{\mathbf{X}}
\def\cN{\mathcal{N}}

\begin{document}
\title{Discovery of Governing Equations with Recursive Deep Neural Networks}
\author[ J. Zhao \& J. Mau ]{ Jia Zhao\affil{1}\comma \corrauth and Jarrod Mau \affil{1}}
\address{ \affilnum{1}\ Department of Mathematics \& Statistics, Utah State University, Logan, UT, USA }
\emails{ {\tt jia.zhao@usu.edu.} (J.~Zhao)}

\begin{abstract}
Model discovery based on existing data has been one of the major focuses of mathematical modelers for decades. Despite tremendous achievements of model identification from adequate data, how to unravel the models from limited data is less resolved.  In this paper, we focus on the model discovery problem when the data is not efficiently sampled in time. This is common due to limited experimental accessibility and labor/resource constraints. 
Specifically, we introduce a recursive deep neural network (RDNN) for data-driven model discovery. This recursive approach can retrieve the governing equation in a simple and efficient manner, and it can significantly improve the approximation accuracy by increasing the recursive stages. In particular, our proposed approach shows superior power when the existing data are sampled with a large time lag, from which the traditional approach might not be able to recover the model well. Several widely used examples of dynamical systems are used to benchmark this newly proposed recursive approach. Numerical comparisons confirm the effectiveness of this recursive neural network for model discovery. 
\end{abstract}

\keywords{Deep Neural Networks; Physics Informed Neural Networks; Model Discovery; Dynamical Systems}
\maketitle

\section{Introduction}
With the recent advancements of computing power and integrated machine learning software development, deep neural network (DNN) has emerged as a powerful tool in the research communities and industrial applications. By leveraging a large amount of available data, the deep neural network has been widely appreciated to tackle many data-driven scientific problems that are otherwise not accessible, significantly advancing the knowledge in some fields that are less unraveled before. One area that draws particular attention from the applied mathematicians is utilizing deep neural networks to solve differential equations, especially high-dimensional problems. And noticeably, many advancements are achieved so far, to name a few \cite{E, deep-energy, discussion-pde, complex-geometries, Adaptive, DGM, HanEPNAS2018, Zhao-Colby}.  In the meanwhile, as an inverse problem, how to discover the hidden governing differential equations from the existing data is the other primary focus, to name a few seminal work \cite{RK-DNN, RaissiJCP, DeepXDE, DHPM, Xiu-JCP, Multiscale, multistep, OnsagerNet, Symplectic,PDE-Net,PDE-discovery}. Some other related focus is on how to define the neural network structure and how to design strategies to train them properly. In this paper, we mainly focus on model discovery from data with deep neural networks.

Mathematically speaking, a deep neural network is a mapping, with specified input and output structure/dimensions. Here the mapping is defined as a composition of nonlinear functions, with a large number of free parameters in the functions.  This mapping ideally can approximate any function to any order of accuracy, so long as the freedom is large enough, which is also known as universal approximation theory\cite{ApproximationTheory}. In the past few years, many seminal works have been published, focusing on model identification from data with machine learning, in particular deep neural networks. Roughly, these works can be categorized into several types. First of all, or maybe the simplest case, one can assume the model is a linear combination of specific terms. The undetermined coefficients can be recovered by the least square approximations based on the available data. For instance, the authors in \cite{RaissiJCP} recover the undetermined coefficients in the PDE models by the physics informed neural networks(PINNs). And in \cite{Brunton3932}, the authors propose a dictionary for all possible terms and then identify the coefficients through a sparse (lasso) regression.  In the second case, one may represent a flow map by a deep neural network. Once the flow map is identified by training the neural network, the dynamic can be advanced in time by the flow map.  For instance, in  \cite{Xiu-JCP}, the authors propose a deep neural network to approximate the flow map with fixed time lag for dynamical systems, such that the dynamics can be predicted by using the flow map recursively.  Some improvements are introduced in followed up papers, such as treating non-autonomous systems \cite{Xiu-Nonautonomous}, defining a residual neural network based on existing coarse-grain approximation of the flow map \cite{Xiu-Residual}, and dealing with multi-scale problems with multiple flow maps \cite{Multiscale}.  In the third case, one may approximate the right-hand side of the system, particularly the reactive kinetics of the dynamically systems, by deep neural networks. Once the reactive kinetics are identified, the dynamical systems can be solved by treating the reactive kinetics as black-box functional. For instance, the authors in \cite{RK-DNN} introduce the Runge-Kutta neural network to approximate the dynamics systems, and the authors in \cite{multistep} introduce the multistep methods for dynamical system approximation. In \cite{Zhao-Image}, we introduced the pseudo-spectral physics informed neural networks to identify the bulk term in the phase-field equations.

However, many of the existing approaches have relatively strong assumptions/requirements on the data, which sometimes are not practical in reality.  In this paper, we introduce a recursive deep neural network with weaker requirements on the existing data, while retrieving the governing systems accurately. Specifically, unlike the requirements in \cite{Xiu-JCP, Xiu-Residual}, we don't necessarily need the data to be sampled with uniform time lag; unlike the requirements in \cite{multistep}, our method works for data from different time senescences, and unlike the requirements in \cite{RK-DNN}, out method works for data that are sampling with big time lags. It turns out our proposed recursive deep neural network works for more scenarios with weaker restrictions on the available data, making it more closely related to the data collected in the real world.  

The rest of this paper is organized as follows. In Section 2, we will explain the problem set up and introduce the recursive deep neural networks for the model discovery. In Section 3, several numerical examples are given to illustrate the effectiveness of the proposed recursive approach. In the last section, a brief discussion and conclusions are provided.

\section{Problem Setups}

\subsection{Statement of the problem}
First of all, we formulate our data-driven model discovery problem in a mathematical manner. Without loss of generality, we consider the unknown parametrized dynamical system
\beq \label{eq:ODE}
\Phi_t = F(\Phi, \mu, t), 
\eeq 
with $\Phi$ the state variables, $F$ the right hand side functional and $\mu$ are parameters. Assume some data measurement of trajectories are accessible, denoted by
\beq
\left\{ \Phi(t_j^{(i)}), \quad  \mu^{(i)},   \quad j=1,2,\cdots, K^{(i)} , \quad i=1,2,\cdots, N_T  \right\},
\eeq 
where $N_T$ is the number of measured time series, and $K^{(i)}$ is the number data entries in each of time series. And the goal is to recover the governing equation of \eqref{eq:ODE} from the data using deep neural networks.

Without loss of generality, we can always re-organize the data as the following data pairs
\beq \label{eq:Data}
\left\{ (\Phi_j^1,\,\,\,  t_j^1, \,\,\, \Phi_j^2,\,\,\,   t_j^2, \,\,\, \mu_j) \right\}_{j=1}^N.
\eeq 
Here $\Phi_j^1, \Phi_j^2$ are snapshots of the solution for \eqref{eq:ODE} with a specific initial condition and parameters $\mu=\mu_j$, at time $t_j^1$ and $t_j^2$ respectively.  And $N$ is the number of data pairs. In particular, unlike the requirements of many existing machine-learning model discovery methods, $\delta_j :=t_j^2 - t_j^1$ is not necessary small, making our approach applicable to a more general case. 
Also, the $N$ data pairs are not necessarily from the same time sequence, weaker than the requirements of many existing approaches as well. So, the problem is phrased as
\beq 
\mbox{Find the dynamical system } \Phi_t = F(\Phi, \mu, t),  \mbox{ given the data } \left\{ (\Phi_j^1, \,\,\,t_j^1, \,\,\, \Phi_j^2, \,\,\,t_j^2,  \,\,\, \mu_j) \right\}_{j=1}^N.
\eeq 

It turns out the problem above can be simplified. We point out the fact that  the parametrized problem in \eqref{eq:ODE} can also be formulated  as
\beq
\Psi_t = \mathcal{F}(\Psi, t), \mbox{ with } \Psi = \left[\begin{array}{l} \Phi \\ \mu \end{array} \right], \quad  \mbox{ and } \mathcal{F} = \left[\begin{array}{l} F(\Psi, t) \\ \mathbf{0} \end{array} \right],
\eeq 
i.e. the parametrized problem can be rewritten as a non-parameterized problem, by introducing extra state variables.
Therefore, we restate the general problem in this paper as:
\beq \label{eq:ODE-2}
\mbox{Find the dynamical system } \Phi_t = F(\Phi, t),  \mbox{ given the data } \left\{ (\Phi_j^1, \,\,\,t_j^1, \,\,\, \Phi_j^2, \,\,\,t_j^2) \right\}_{j=1}^N.
\eeq 

In the rest of this paper, we mainly focus on the problem in \eqref{eq:ODE-2} and propose a novel recursive deep neural network strategy to unravel the dynamical system \eqref{eq:ODE} based on the data pairs \eqref{eq:Data}.

\subsection{Deep Neural Networks}
Mathematically, the feed-forward neural network could be defined as compositions of nonlinear functions. Given an input $x \in \hR^{n_1}$, with $n_1$ the dimension of the input, we can define the general feed-forward deep neural network as \cite{DeepLearningMath} 
\beq \label{eq:NN}
\bea{l}
a^{[1]} = x \in \hR^{n_1}, \\
a^{[k]} = \sigma \Big( W^{[k]} a^{[k-1]} + b^{[k]} \Big) \in \hR^{n_k}, \quad \mbox{ for } k=2,3,\cdots, L,
\eea 
\eeq 
where  $W^{[k]} \in \hR^{n_k \times n_{k-1}}$ and $b^{[k]} \in \hR^{n_k}$ denote the weights and biases at layer $k$ respectively, $\sigma$ denotes the activation function, which is usually nonlinear. Here  $a^{[l]} \in \hR^{n_l}$ is the output of the $l$-th layer, which is the input for $l+1$-th layer.
Essentially neural networks are non-linear mappings with many parameters.  

In this paper, to solve the problem in \eqref{eq:ODE-2}, we introduce a physics-uniformed  feed-forward neural network 
\beq
\cN_F: (\Phi, t;\theta ) \rightarrow \cN_F(\Phi, t; \theta)
\eeq
to approximate the kinetic term $F$ in \eqref{eq:ODE-2}, where $\theta = \{ W^{[k]}, b^{[k]}, k=2,3,\cdots, L \}$ are the free parameters, i.e. the weights and biases  in $\cN_F$ to be tuned. In the rest discussion, we may omit the notation $\theta$ in $\cN_F$ without confusion.

\subsection{Physics informed neural networks}
In the meanwhile, we introduce a physics-informed residual neural network 
\beq
\cN_r: \quad (\Phi_j^1, \,\,\,t_j^1, \,\,\, \Phi_j^2, \,\,\,t_j^2) \rightarrow  \cN_r(\Phi_j^1, \,\,\,t_j^1, \,\,\, \Phi_j^2, \,\,\,t_j^2),
\eeq 
which shares the same parameters as $\cN_f$.
Next, we illustrate how the physically informed neural network $\cN_r$ is defined.

Presumably, if $\cN_F$ is equivalent to $f$, we shall have 
\beq \label{eq:time-integrator}
\Phi_j^2 =\Phi_j^1 + \int_{t_j^1}^{t_j^2} \cN_F(\Phi, t) dt,
\eeq 
from the data $(\Phi_j^1, \,\,\,t_j^1, \,\,\, \Phi_j^2, \,\,\,t_j^2)$.
Unfortunately, there is no analytical formula for the integral, so numerical approximation shall be carried out starting here.   For simplicity,we denote $h_j:=t_j^2 - t_j^1$ in the rest discussion.

Learned from classical numerical methods for differential equations, if we approximate the integral by rectangular rule, we will have 
$$
\Phi_j^2 \approx \Phi_j^1 + h_j \cN_F(\Phi_j^1, t_j^1) + O(h_j^2),
$$
from which we can define $\cN_r$ as
\beq \label{eq:Euler}
\cN_r:  (\Phi_j^1, \,\,\,t_j^1, \,\,\, \Phi_j^2, \,\,\,t_j^2) \rightarrow  \Phi_j^2 - \Phi_j^1 - h_j \cN_F(\Phi_j^1, t_j^1).
\eeq
Similarly, we can define the residual neural network as
\beq
\cN_r  (\Phi_j^1, \,\,\,t_j^1, \,\,\, \Phi_j^2, \,\,\,t_j^2)  \rightarrow  \Phi_j^2 - \Phi_j^1 - h_j \cN_F(\Phi_j^2, t_j^2).
\eeq 
if we use the following rectangular rule
$$\Phi_j^2  \approx  \Phi_j^1 + h_j \cN_F(\Phi_j^2, t_j^2) + O(h_j^2).$$ 
Or if we approximate the integral by a trapezoid rule, we will have 
$$
\Phi_j^2  \approx  \Phi_j^1 + \frac{h_j}{2} (\cN_F(\Phi_j^1, t_j^1) + \cN_F(\Phi_j^2, t_j^2) ) + O(h_j^3).
$$
Then we can define $\cN_r$ as
\beq\label{eq:CN}
\cN_r:  (\Phi_j^1, \,\,\,t_j^1, \,\,\, \Phi_j^2, \,\,\,t_j^2)  \rightarrow \Phi_j^2 - \Phi_j^1 - \frac{h_j }{2} (\cN_F(\Phi_j^1, t_j^1) + \cN_F(\Phi_j^2, t_j^2) ).
\eeq 

Once the two neural networks: $\cN_F$ and $\cN_r$ are defined, the problem in \eqref{eq:ODE-2} can be solved by tackling the optimization problem 
\beq \label{eq:optimization}
\mbox{ Find } \cN_F(\Phi, t; \theta) \mbox{ such that } \min_{\theta \in \cN_F} \sum_{j=1}^N \| \cN_r (\Phi_j^1, \,\,\,t_j^1, \,\,\, \Phi_j^2, \,\,\,t_j^2) \|^2,
\eeq 
where $\cN_F$ will be the approximation for $f$.

In principle, if $\max_{j} |t_j^2-t_j^1|$ is small enough (near zero), the reformulated optimization problem \eqref{eq:optimization} will be consistent with \eqref{eq:ODE-2}, and thus $\cN_f$ provides an accurate approximation for the kinetic term $f$ in \eqref{eq:ODE-2}. But unfortunately, this is not practical in reality, as $|t_j^2-t_j^1|$ is usually large in the data, saying due to the sampling capability in experiments or resource constraints. When $|t_j^2-t_j^1|$ is large, the truncation errors in the integral approximation, such as in \eqref{eq:Euler} or \eqref{eq:CN}, are not negligible any more. Thus, a better definition for the residual neural network $\cN_r$ is necessary.

\subsection{Physics-informed Recursive Neural Networks}
In this sub-section, we introduce the physics-informed recursive neural network to define $\cN_r$. Instead of approximating the integral in \eqref{eq:time-integrator} directly, we approximate it in a recursive approach. Without loss of generality, we discretize the time interval $[t_j^1, t_j^2]$ into small segments, $t_j^1=\tau_0 < \tau_1 < \tau_2 < \cdots < \tau_M=t_j^2$, $j=1,2,\cdots N$,  with $M$ the number of  segments (as a hyper-parameter).

As a recursive version of \eqref{eq:Euler}, the residual neural network can be defined as
\beq \label{eq:Recursive-Euler}
\cN_r: (\Phi_j^1, \,\,\,t_j^1, \,\,\, \Phi_j^2, \,\,\,t_j^2) \rightarrow  \Phi_j^2 - \Phi_{j, M}
\eeq
where $\Phi_{j,M}$ is defined by the following recursive formula
\beq
\bea{l}
\Phi_{j,0} = \Phi_j^1, \\
\Phi_{j, s+1} = \Phi_{j, s}  + (\tau_{s+1}- \tau_s)  \cN_f(\Phi_{j, s}, \tau_s),\quad  s=0,1,\cdots, M-1.
\eea
\eeq 
And it can be easily shown the local truncation error (LTE) for this approximation is $O((\frac{t_j^2-t_j^1}{M})^2)$.
Then the optimization problem in \eqref{eq:optimization} can be solved with this new $\cN_r$ defined  in \eqref{eq:Recursive-Euler}.

To further reduce the local truncation error, one can introduce higher-order approximation for the integral. For instance, inspired by the classical 4-stage Runge-Kutta method, we can introduce  the recursive Runge-Kutta neural network
\beq \label{eq:Recursive-RK}
\cN_r: (\Phi_j^1, \,\,\,t_j^1, \,\,\, \Phi_j^2, \,\,\,t_j^2) \rightarrow  \Phi_j^2 - \Phi_{j, M}.
\eeq
where $\Phi_{j,M}$ is defined  by the recursive formula 
\beq 
\bea{l}
\Phi_{j,0} = \Phi_j^1, \\
\Phi_{j, s+1} = \Phi_{j, s} + \frac{h_j}{6} (K_{s1} + 2K_{s2} + 2K_{s3} +K_{s4}), \quad  h_j = \tau_{s+1}-\tau_s, \\
K_{s1} = \cN_F(\Phi_{j,s}, \tau_s )， \\
K_{s2} = \cN_F(\Phi_{j,s} +\frac{h_j}{2}K_{s1} , \tau_s + \frac{h_j}{2})， \\
K_{s3} = \cN_F(\Phi_{j,s} +\frac{h_j}{2}K_{s2} , \tau_s + \frac{h_j}{2})， \\
K_{s4} = \cN_F(\Phi_{j,s} + h_j K_{s3}, \tau_{s+1}),\quad s=0,1,\cdots, M-1.
\eea
\eeq 
The four-stage RK method is a fourth-order method, meaning that the local truncation error for this recursive approach is $O((\frac{t_j^2-t_j^1}{M})^5)$. Similarly, the optimization problem in \eqref{eq:optimization} can be solved with the residual network defined in \eqref{eq:Recursive-RK}.

We remark that this idea is not limited to the two approaches above. Here, we present a pseudo code to define this recursive physics-informed residual neural network, as shown below.
\begin{center}
\begin{lstlisting}[language=Python, caption=Definition of Recursive Physics-informed Residual Neural Networks]
def neural_net(self, X, weights, biases):
    num_layers = len(weights) + 1
    H = X
    for l in range(0, num_layers - 2):
        W = weights[l]
        b = biases[l]
        H = tf.tanh(tf.add(tf.matmul(H, W), b))
    W = weights[-1]
    b = biases[-1]
    Y = tf.add(tf.matmul(H, W), b)
    return Y

def net_f(self, u, t):
    return self.neural_net(tf.concat([u,t],1), self.weights, self.biases)
	
def time_march(phi_s1, tau1, tau2, method):
	h = tau2 - tau1

	if method == "Euler":
		return phi_s1 + h * net_f(phi_s1, tau1)

	elif method == "RK2":
		K1 = net_f(phi_s1, tau1)
		K2 = net_f(phi_s1 + h/2 * K1, tau1 + h/2)
		return phi_s1 + h * K2

	elif method == "RK3":
		K1 = net_f(phi_s1, tau1)
		K2 = net_f(phi_s1 + h * K1, tau2)
		K3 = net_f(phi_s1 + h/4* (K1+K2), tau1 + h/2)
		return phi_s1 + h/6 * (K1 + K2 + 2*K3)

	elif method == "RK4":
		K1 = net_f(phi_s1, tau1)
		K2 = net_f(phi_s1 + h/2 * K1, tau1 + h/2)
		K3 = net_f(phi_s1 + h/2 * K2, tau1 + h/2)
		K4 = net_f(phi_s1 + h * K3 + tau2)
		return phi_s1 + h/6 * (K1 + 2*K2 +2*K3 + K4)
	else method == "RK5":
        K1 = net_f(phi_s1, tau1)
        K2 = net_f(phi_s1 + 1/4 * K1, tau1 + h/4)
        K3 = net_f(phi_s1 + 1/8 * K1 + 1/8 * K2, tau1 + h/4)
        K4 = net_f(phi_s1 - 1/2 * K2 + K3, tau1 + h/2)
        K5 = net_f(phi_s1 + 3/16 * K1 + 9/16 * K4, tau1 + h*3/4)
        K6 = net_f(phi_s1-3/7*K1+2/7*K2+12/7*K3-12/7*K4+8/7*K5, tau2)
        return phi_s1+h/90*(7*K1 + 32*K3 + 12*K4 + 32*K5 + 7*K6)
	else:
		return None

def residual_net(phi1, phi2, t1, t2, method):
	phi_s1 = phi1
	... # define time steps: tau based on t1,t2.
	for i in range(M):
		phi_s2 = time_march(phi_s1, tau[i], tau[i+1], method)
		phi_s1 = phi_s2
	return phi2 - phi_s2
\end{lstlisting}
\end{center}

\section{Numerical Examples}
Next, we present several numerical examples to illustrate the effectiveness of the proposed recursive approach on uncovering dynamical systems from data. 
Note the theory presented in this paper is rather general, and it applies to a variety of non-autonomous, parametric dynamical systems, and the data is allowed to be sampled from non-uniform time intervals. For simplicity of illustration purpose, we only consider autonomous ODE systems, and we sample the data with a uniform time step size $\Delta t$ and large time lags. More complicated scenarios will be further investigated in our later research. 

The numerical tests in this section are conducted in the following steps:
\begin{itemize}
\item As proof of principle, we generate synthetic data by solving some classical dynamic systems, which we name the true models.  The generated data set is provided in the form of \eqref{eq:ODE-2}. Mainly, we fix a sampling domain $\Omega$, and randomly sampled for the initial conditions $\{\Phi_i^1\}_{i=1}^N \subset \Omega$, and use the true model to advance the solution by a time lag $\Delta t$, to obtain $\{\Phi_i^2\}_{i=1}^N$, i.e., the data are in the format of
$\{ (\Phi_i^1, \Phi_i^2, \Delta t)  \}_{i=1}^N$. 
\item  Next, we construct the physics uninformed neural network $\cN_F$ and the physics-informed recursive neural network $\cN_r$. Then we solve the optimization problem \eqref{eq:optimization} for model discovery from the synthetic data above, without using the knowledge of the true models. 
\item Once the neural networks are well trained, the neural network $\cN_F$ would be a good approximation to $f$ in \eqref{eq:ODE-2}. Given a random initial condition, the predictions of the dynamic using the learned model  $\Phi_t = \cN_F(\Phi, t)$ and the true model in \eqref{eq:ODE-2} will be compared. 
\end{itemize}

In all the numerical examples, we fixed some hyper-parameters.
For the neural network construction of $\cN_F$, we use a feed-forward neural network with one hidden layer of 128 neurons and use $\tanh$ as the activation function. For the recursive neural network $\cN_r$, we use the one defined in \eqref{eq:Recursive-RK}.  In all the simulations below, otherwise specified, we will use the Adam method for 10000 steps of optimization, followed directly by an L-BFGS-B method to fine-tune the result.  Notice a fine-tuning of the hyper-parameters, such as the neural network structure (saying the width and depth), the training rate, the mini-batches, weights on the loss function, etc., will improve the accuracy, but these techniques are out of our current research interest. Interested readers are highly encouraged to explore.

\subsection{Two-dimensional damped harmonic oscillator}
First of all, we start with a simple example. Consider the two-dimensional damped harmonic oscillator with cubic reactive dynamics, which is given as
\beq \label{eq:cubic}
\dot{\bX} = \left[
\begin{array}{ll}
-0.1 & \,\,\,\, 2.0 \\
-2.0 & -0.1
\end{array} 
\right] 
\bX^3,  \quad \mbox{ where } \bX=\left[\begin{array}{l} x \\ y  \end{array}\right].
\eeq 

In this example, we pick the sampling domain $\Omega = [-2.5 ,\,\,\, 2.5] \times [-2.5, \,\,\, 2.5]$, and sample  $N=1000$ points from $\Omega$ as $\{\Phi_j^1 \}_{j=1}^N$ using Latin hypercube sampling,  fix the corresponding time lag $\Delta t$, and calculate $\{\Phi_j^2 \}_{j=1}^N$ using the true model. Once $\cN_F$ is trained with the data, we obtain the learned model $\bX_t = \cN_F(\b X)$. we randomly choose an initial condition. Here we pick $\bX= [2,  \,\,\, 0]^T$. Then we use both models (the discovered model and the true model) to calculate the dynamics for $t \in [0, 25]$. 

The relative $l^2$ errors between accurate solutions solved using the true model and predicted solutions solved using the learned models, which are learned from data with different time lag size $\Delta t$ and recursive stages $M$, are summarized in Table \ref{tab:Cubic}. We observe that when the time lag size for the training data in \eqref{eq:ODE-2} is small, even a single recursive stage ($M=1$) works fine. However, when the time lag size is large in the training data, such as $\Delta t = 0.2$, a neural network with more recursive stages, i.e. larger M, provides more accurate approximation. In other words, our proposed recursive deep neural network shows superior approximation power than the classical approaches when the available training data are sampled with a large time lag.
\begin{table}[H]
\caption{relative $l^2$ errors for two-dimensional damped harmonic oscillator in \eqref{eq:cubic} calculated with the learned model $\dot{\bX}=\cN_F(\bX)$.}
\begin{center}
\begin{tabular}{ |c |c |c |c|c|}
\hline
& M=1 & M=2  & M = 5 & M = 10\\ 
\hline
$\Delta t$ = 0.01&0.0135 & 0.0025& 0.0162 & 0.0108 \\  
\hline
$\Delta t$ = 0.05& 0.0201& 0.007 &   0.0025 &  0.0031  \\
\hline
$\Delta t$ = 0.1 & 0.21663 & 0.0608 &  0.0081 & 0.0039   \\
\hline
$\Delta t$ = 0.2 &0.9753&0.0879 &  0.0059 &  0.0060 \\
\hline
\end{tabular}
\label{tab:Cubic}
\end{center}
\end{table}

\begin{figure}[H]
\centering
\subfigure[$\Delta t=0.2$, $M=1$]{\includegraphics[width=0.475\textwidth]{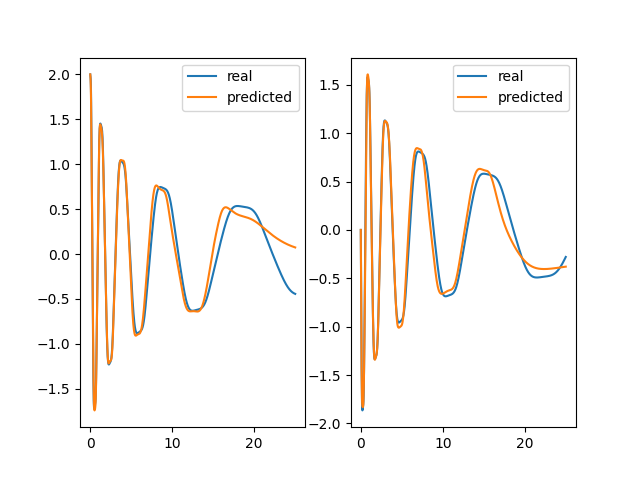}}
\subfigure[$\Delta t=0.2$, $M=5$]{\includegraphics[width=0.475\textwidth]{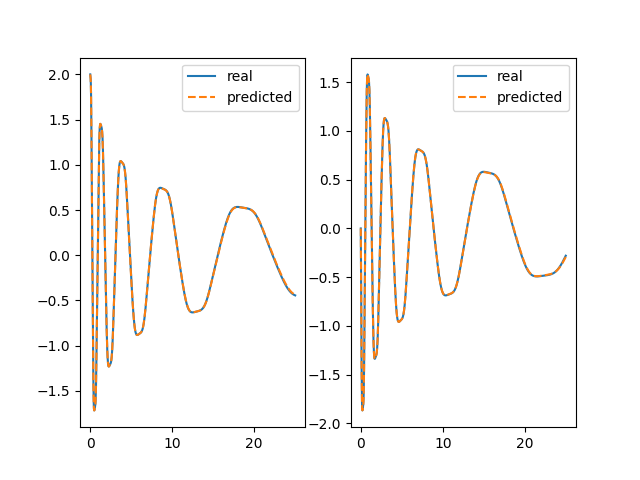}}
\caption{Numerical predictions of the two-dimensional damped oscillator in \eqref{eq:cubic} with the learned dynamical system $\dot{\bX} = \cN_F(\bX)$. In (a), prediction by the learned model using a  single stage ($M=1$)  neural network. In (b) prediction by the learned model using a $M=5$ stage neural network.Here the orange lines represent the predicted solutions by the learned model, and the blue lines indicate the accurate solutions that are obtained by solving the true model.}
\label{fig:cubic}
\end{figure}

Furthermore, we also plot the predicted solutions using the learned models $\dot{\bX}=\cN_F(\bX)$ of various recursive stages that are trained by data with a time lag $\Delta t=0.2$, as shown in Figure \ref{fig:cubic}. We observe the learned model using $M=5$ recursive stage deep neural network can accurately capture the dynamics, while the classical single-stage approach fails. We remark that this system has also been studied in \cite{multistep}, where the multistep neural network can capture the dynamics with $\Delta t = 0.01$. Here our recursive approach can capture the dynamics well even with $\Delta t = 0.2$.

\subsection{Glycolytic oscillator}
Next, we consider a more complicated nonlinear dynamic system from \cite{Glycolytic}, which has also been investigated in \cite{Multiscale}. This model consists of seven coupled ODEs, for the concentration of 7 biochemical species \cite{Glycolytic}, given by
\beq \label{eq:glycolytic}
\bea{l}
\frac{d S_1}{dt} = J_0 - \frac{k_1 S_1 S_6}{1 + (S_6/K_1)^q}, \\
\frac{d S_2}{dt} = 2 \frac{k_1 S_1 S_6}{1+(S_6/K_1)^q} - k_2 S_2 (N-S_5) - k_6 S_2 S_5, \\
\frac{d S_4}{dt} = k_3 S_3(A-S_6) - k_4 S_4 S_5 - \kappa(S_4 - S_7), \\
\frac{dS_5}{dt} = k_2 S_2 (N-S_5) - k_4 S_4 S_5  - k_6 S_2 S_5, \\
\frac{d S_6}{dt} = -2 \frac{k_1 S_1 S_6}{1 + (S_6/K_1)^q} + 2k_3 S_3 (A-S_6) - k_5 S_6, \\
\frac{ d S_7}{dt} = \psi_k (S_4 - S_7) -k S_7.
\eea 
\eeq 
with the coefficient chosen directly from  from \cite{Glycolytic}, given as $J_0=2.5$, $k_1 =100$, $k_2=6$, $k_3=16$, $k_4=100$, $k_5 = 1.28$, $k_6 = 12$, $k=1.8$, $\kappa=13$, $q=4$, $K_1 = 0.52$, $\psi=0.1$, $N=1$, and $A=4$.
Note \eqref{eq:glycolytic} is an oscillatory system.

We choose the sampling domain $\Omega=[0, 2] \times[0, 3] \times [0, 0.5 ] \times [0, 0.5] \times [0, 0.5] \times [0.14, 2.67] \times [0.05, 0.15]$, and  using Latin hypercube sampling to generate 1000 points as $\{\Phi_i^1\}_{i=1}^N$, and calculate $\{\Phi_i^2\}_{i=1}^N$ with the true model accordingly. After the neural network $\cN_F$ is trained, we obtained the learned model $\Phi_t = \cN_f(\Phi)$. 
We pick a random initial condition $[1.1, 1.0, 0.075, 0.175, 0.25, 0.9, 0.095]$, the same as \cite{multistep}, to calculate solutions in the time interval $[0, 5]$. 

The $l^2$ errors between  the accurate solutions and the predictions using learned models that are trained from scenarios with different time lag steps and recursive stages are summarized in Table \ref{tab:glycolytic}. A similar phenomenon as the previous example is observed. Mainly, when the time lag for the data is large, a single-stage neural network leads to large errors. In the meanwhile, a recursive neural network with more intermediate stages provides a reliable approximation. 
\begin{table}[H]
\caption{$l^2$ error of the solutions predicted by the learned model  for the glycolytic oscillator systems in \eqref{eq:glycolytic} }
\begin{center}
\begin{tabular}{ |c |c |c |c|c|}
\hline
& M=1 & M=2  & M = 5 & M = 10\\ 
\hline
$\Delta t$ = 0.2&0.7105 & 0.1418&  0.0281&  0.0034 \\  
\hline
$\Delta t$ = 0.5&  1.3040 &  0.1090&   0.1976 &  0.0366  \\
\hline
\end{tabular}
\end{center}
\label{tab:glycolytic}
\end{table}

\begin{figure}
\centering
\subfigure[$\Delta t=0.2$, $M=1$]{\includegraphics[width=0.45\textwidth]{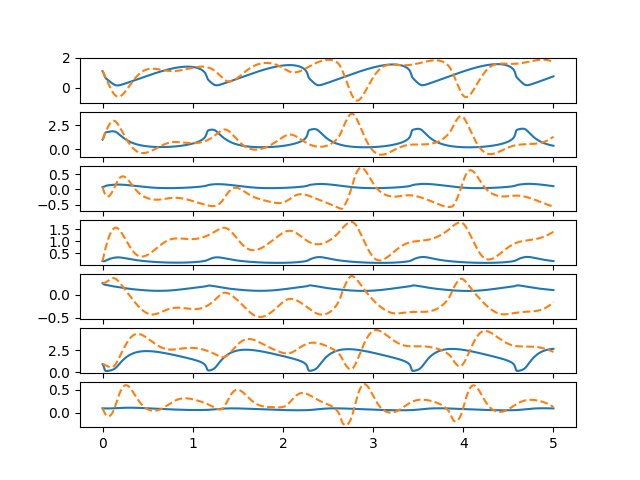}}
\subfigure[$\Delta t=0.2$, $M=2$]{\includegraphics[width=0.45\textwidth]{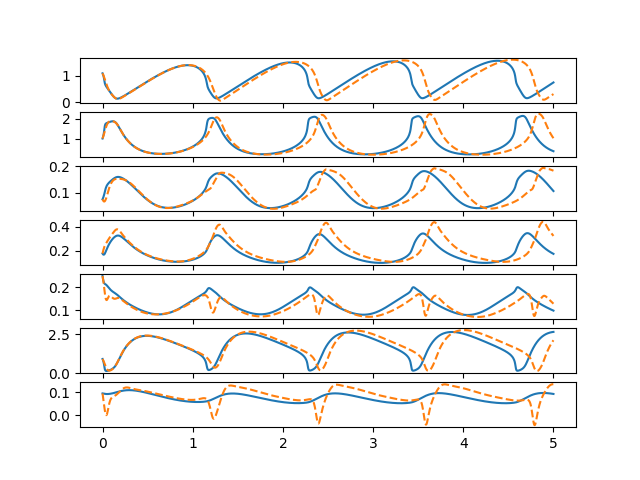}}

\subfigure[$\Delta t=0.2$, $M=5$]{\includegraphics[width=0.45\textwidth]{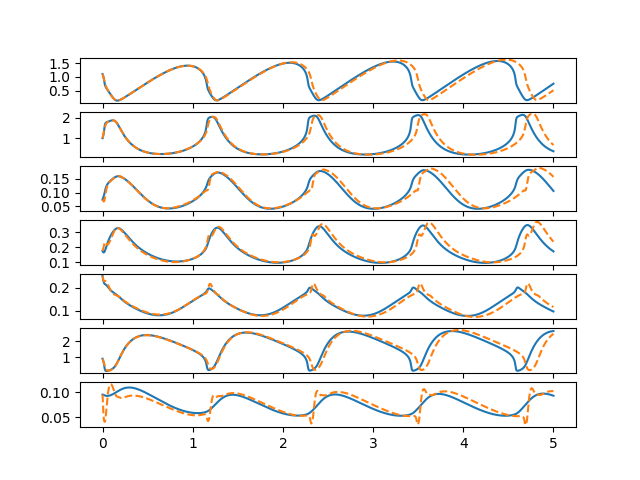}}
\subfigure[$\Delta t=0.2$, $M=10$]{\includegraphics[width=0.45\textwidth]{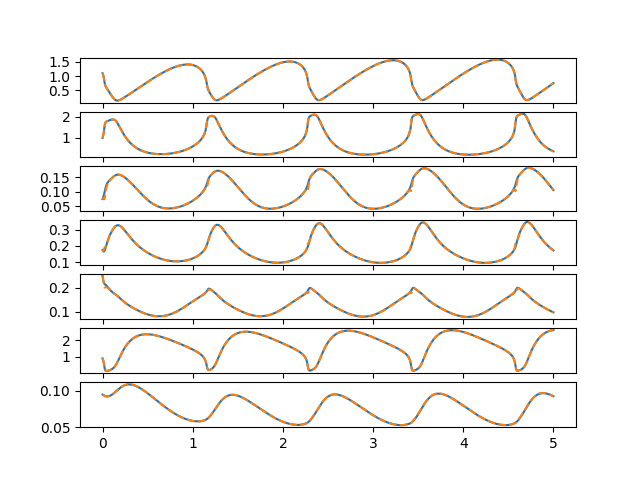}}
\caption{Numerical predictions of nonlinear complex biological systems of \eqref{eq:glycolytic} with the learned model. Here the model is trained with data that has time lag of $\Delta t =0.2$. The models of different intermediate stage $M$ are compared. And we observe that the model with big intermediate stage $M$ can predict the dynamics well. Here the orange lines represent the predicted solutions by the learned model, and the blue lines indicate the accurate solutions that are obtained by solving the true model} 
\label{fig:glycolytic}
\end{figure}

In addition, the numerical predictions using various trained models with data lag $\Delta t = 0.2$ are summarized in Figure \ref{fig:glycolytic}. We observe the learned models by using the deep neural network with more recursive stages, such as $M=5$, provide reliable predictions for the oscillation dynamics.

\subsection{Hopf bifurcation}

So far, we have illustrated that the learned models by the recursive neural network approach can provide an accurate approximation for the kinetics of the nonlinear and complicated autonomous dynamical systems. Thus this approach can discover the dynamical systems through data successfully. To further push the application limit of our proposed approach, we investigate a parametric dynamical system in this example. Consider the following parametrized dynamical system, which has also been used in \cite{multistep}, as
\beq \label{eq:Hopf}
\bea{l}
\dot{x} = \mu x + y - x(x^2+y^2), \\
\dot{y} = -x + \mu y - y (x^2 + y^2),
\eea 
\eeq 
where $\mu$ is a parameter. It is known this system goes through bifurcation as the parameter $\mu$ varies. In particular, there is a transition from the fixed point for $\mu < 0$ to the limit cycle for $\mu > 0$. 

We point out that this parameterized system can be equivalently formulated as
\beq
\left[
\begin{array}{l}
\mu \\
x \\
y
\end{array}
\right]_t =
\left[
\begin{array}{l}
0 \\
\mu x + y - x(x^2+y^2) \\
-x + \mu y - y (x^2 + y^2)
\end{array}
\right].
\eeq 
So, in this context, $\Phi=(\mu \,\,\, x \,\,\, y)^T$ for the form in \eqref{eq:ODE-2}.

In this example, we consider the domain $\Omega=[-1, \,\,\, 1] \times [-2, \,\,\, 2] \times [-1, \,\,\, 1]$, and use Latin hypercube sampling to randomly sample 1000 points from $\Omega$ as $\{ \Phi_i^1\}_{i=1}^N$, and calculate the corresponding $\{\Phi_i^2\}_{i=1}^N$ with the true model. After the model is trained, we choose the same initial conditions as \cite{multistep}, and calculate the solutions for $t\in [0, 75]$. 

The $l^2$ errors using different trained models are summarized in Table \ref{tab:Hopf}.  We observe that when the time step is large, such as $\Delta t = 2$, the single-stage model fails to predict accurate dynamics, but the recursive model can predict reliable dynamics.  
\begin{table}[H]
\caption{$l^2$ error of the predicted solutions using the learned model with different settings for the parametrized model in \eqref{eq:Hopf}}
\begin{center}
\begin{tabular}{ |c |c |c |c|c|}
\hline
& M=1 & M=2  & M = 5 & M = 10\\ 
\hline
$\Delta t$ = 0.5 & 0.0875 &  0.0312&  0.0257 &  0.007 \\  
\hline
$\Delta t$ = 1.0 &  0.2439&  0.0469&   0.0076 &  0.0114  \\
\hline
$\Delta t$ = 2.0 & 1.0732 & 0.2287 & 0.0321 & 0.0129 \\
\hline
\end{tabular}
\label{tab:Hopf}
\end{center}
\end{table}

\begin{figure}[H]
\centering
\subfigure[$\Delta t=2$, $M=1$]{\includegraphics[width=0.475\textwidth]{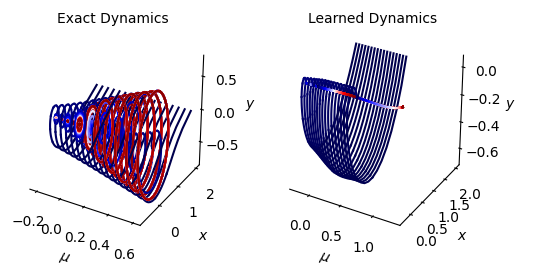}}
\subfigure[$\Delta t=2$, $M=2$]{\includegraphics[width=0.475\textwidth]{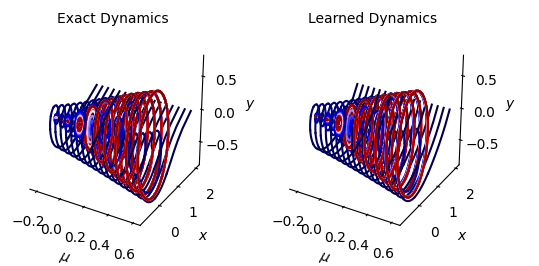}}

\subfigure[$\Delta t=2$, $M=5$]{\includegraphics[width=0.475\textwidth]{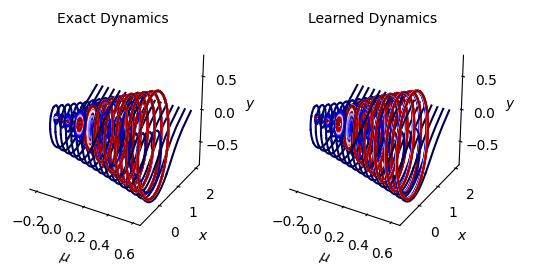}}
\subfigure[$\Delta t=2$, $M=10$]{\includegraphics[width=0.475\textwidth]{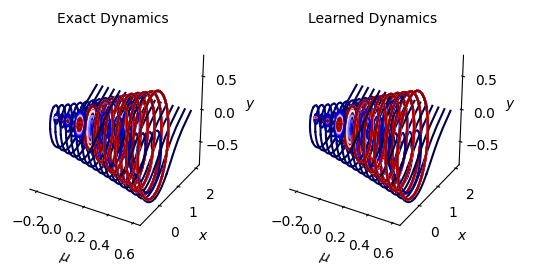}}
\caption{Numerical predictions of Hopf bifurcation using the trained models with various recursive stages. This figure shows that even when the training data are collected with time lag $\Delta t=2$, the recursive approach, such as $M=5, 10$, can accurately capture the bifurcation transitions.}
\label{fig:Hopf}
\end{figure}

In particular, some predictions using trained models with $\Delta t = 2$ and different stages are summarized in Figure \ref{fig:Hopf}. We remark the, as shown in \cite{multistep}, the multistep method can roughly handle $\Delta t = 0.1$, i.e., capture the transition trend, but not accurately. But our approach captures the transition dynamics accurately, even with a time step $\Delta t = 2$, as shown in Figure \ref{fig:Hopf}(d). This is a large advancement.

\section{Conclusion}
In this paper, we introduce a novel recursive deep neural network for model identification of differential equations. Unlike the idea of flow map discovery for dynamical systems, we present a neural network as a functional to approximate the kinetic part, i.e., the right-hand side of \eqref{eq:ODE}.  Inspired by the numerical methods for differential equations, we formulate the residual neural network recursively, while in each step, there is a multi-stage neural network mapping. 

The proposed recursive neural network has several advantages, making it stand out. First of all, it is easy to implement, and it turns out to improve the accuracy of approximation by simply adjust the intermediate stages. Secondly, it has weaker requirements on the existing data than many current model identification approaches. Thirdly, this idea can be easily combined with many existing model identification strategies to further strengthen the approximation power.  Several numerical examples are presented to illustrate its effectiveness. 
A further discussion on the hyper-parameters and its applications to discover partial differential equations will be investigated in our follow up research.

\section*{Acknowledgments}
The authors would like to acknowledge the support from NSF-DMS-1816783 and NVIDIA Corporation for their donation of a Quadro P6000 GPU for conducting some of the numerical simulations in this paper.

\section*{Conflict of Interest Statement}
On behalf of all authors, the corresponding author states that there is no conflict of interest.

\bibliographystyle{unsrt}

\end{document}